\documentclass[12pt,letterpaper,twoside]{article}

\usepackage{amsmath,amssymb}
\usepackage{amsthm}
\usepackage{stmaryrd}
\usepackage{fancyhdr}
\usepackage{times}
\usepackage{mathrsfs} 
\usepackage{sectsty}
\usepackage[pdftex,dvips]{geometry}

\input diagxy.tex

\sectionfont{\normalsize\bfseries\scshape\centering\MakeUppercase }

\subsectionfont{\normalsize\bfseries\scshape }

\makeatletter
  \def\@seccntformat#1{\csname the#1\endcsname.\quad}
\makeatother

\newtheoremstyle{fact}
     {\topsep}
     {\topsep}
     {\slshape}
     {}
     {\bfseries}
     {}
     { }
     {\thmname{#1}\thmnumber{ #2.}\thmnote{ \rm (#3)}}

\newtheorem{theorem}{Theorem}[section]
\newtheorem*{theorem*}{Theorem} 
\newtheorem{lemma}[theorem]{Lemma}
\newtheorem{proposition}[theorem]{Proposition}
\newtheorem{corollary}[theorem]{Corollary}

\theoremstyle{definition}

\newtheorem*{remark*}{Remark}

\newtheorem*{question*}{Question}
\newtheorem*{examples*}{Examples}  
\newtheorem{example}[theorem]{Example}
\newtheorem*{example*}{Example}
\newtheorem{convention}[theorem]{Convention}
\newtheorem*{convention*}{Convention}

\theoremstyle{fact}

\newtheorem{ftheorem}[theorem]{Theorem}
\newtheorem{flemma}[theorem]{Lemma}

\newtheorem{fcorollary}[theorem]{Corollary}

\def\proofont{\fontseries{bx}\fontshape{sc}\selectfont}
\def\proofname{Proof. }

\newcommand{\rest}{\mbox{\parbox[t]{0.1cm}{$|$\\[-10pt] $|$}}}

\newcommand{\Note}[1]{}

\makeatletter
\renewenvironment{proof}[1][\proofname]{\par
  \normalfont
  \topsep6\p@\@plus6\p@ \trivlist
  \item[\hskip\labelsep\noindent\proofont #1]\ignorespaces
}{%
  \qed\endtrivlist
}
\makeatother

\author{G\'abor Luk\'acs
\thanks{I gratefully acknowledge the generous financial support received
from the Alexander von Humboldt Foundation, the Killam Trusts, and 
Dalhousie University that enabled me to do this research.}}
   
\title{Precompact abelian groups and topological annihilators
\thanks{2000 Mathematics Subject Classification: 22C05 54A20 54D60 (06A15 22A05 54H11)}}

\begin{document}

\makeatletter
\let\mytitle\@title
\chead{\small\itshape G. Luk\'acs / \mytitle }
\fancyhead[RO,LE]{\small \thepage}
\makeatother

\maketitle

\def\thanks#1{} 

\thispagestyle{empty}

\nocite{ComfHernTrig}
\nocite{ComfRem} 
\nocite{HernMaca} 
\nocite{Reid} 
\nocite{Varop}

\begin{abstract}
For a compact Hausdorff abelian group $K$ and its subgroup
$H\leq K$, one defines the
{\em $\mathfrak{g}$-closure} $\mathfrak{g}_K(H)$ of $H$ in $K$
as the subgroup consisting of $\chi \in K$ such that 
$\chi(a_n)\longrightarrow 0$ in
$\mathbb{T}=\mathbb{R}/\mathbb{Z}$ for
every sequence $\{a_n\}$ in $\hat K$
(the Pontryagin dual of $K$) that converges to $0$ in the topology
that $H$ induces on $\hat K$. We prove that every countable subgroup of 
a compact Hausdorff group is $\mathfrak{g}$-closed, and thus give a 
positive answer to two problems of Dikranjan, Milan and Tonolo.
We also show that every $\mathfrak{g}$-closed subgroup
of a compact Hausdorff group is realcompact.
The techniques developed in the paper are used to
construct a close relative of the closure operator $\mathfrak g$ that 
coincides with the $G_\delta$-closure on compact Hausdorff abelian groups, 
and thus captures realcompactness and pseudocompactness of subgroups.
\end{abstract}

\section{Introduction}

\label{sect:intro}

For a sequence $\underline u=(u_n)$ of integers, one defines
$t_{\underline u}(\mathbb T) =\{x \in \mathbb{T} \mid u_n x \longrightarrow 0\}$,
the subgroup of {\em topologically $\underline u$-torsion} elements in
$\mathbb{T}=\mathbb{R}/\mathbb{Z}$, and sets
$\mathfrak{t}(H)=\bigcap\{t_{\underline u}(\mathbb{T}) \mid
H \leq t_{\underline u}(\mathbb{T}) \}$ for every subgroup
$H \leq \mathbb{T}$ (cf. \cite{Diktoptor}, \cite{BarDikMilWeb}, \cite{DikMilTon}).
This enables one to speak of $\mathfrak{t}$-closed subgroups 
(i.e., $H \leq \mathbb{T}$ such that $\mathfrak{t}(H)=H$). The subgroups
$t_{\underline u}(\mathbb{T})$ are called {\em basic} $\mathfrak{t}$-closed
subgroups.
B\'{\i}r\'o, Deshouillers and  S\'os  proved:

\begin{ftheorem}[{\cite[Theorem~2]{BirDesSos}}] \label{fthm:BirDesSos}
Every countable subgroup of $\,\mathbb{T}$ is basic $\mathfrak{t}$-closed.
\end{ftheorem}

Although it is possible to define $t_{\underline u}(G)$ for an arbitrary
abelian topological group in a similar fashion, the closure $\mathfrak{t}_G$ 
thus obtained turned out to have worse properties than anticipated: 
Dikranjan and Di Santo  \cite{DikSan} showed that if $L$ is a non-discrete 
locally compact Hausdorff abelian group and $L \not\cong \mathbb{T}$, 
then there is $x\in L$ such that the subgroup
$\langle x \rangle$ is not even $\mathfrak{t}$-closed in $L$.
In particular, Theorem~\ref{fthm:BirDesSos} fails 
for $L \not\cong \mathbb{T}$.

Following Dikranjan, Milan and Tonolo \cite[2.1]{DikMilTon}, given an abelian 
topological group $G$ and a sequence $\underline{u}=(u_n)$
in its Pontryagin dual $\hat G$, we put 
$s_{\underline u}(G)=\{x\in G\mid u_n(x)\longrightarrow 0
\text{ in $\mathbb{T}$}\}$, the subgroup of {\em topologically
${\underline u}$-annihilating} elements of $G$. For every
subgroup $H \leq G$, we set 
\begin{equation}
\mathfrak g_G(H)=\bigcap\{ s_{\underline u}(G) \mid 
\underline{u}\in\hat G^\mathbb{N},H \leq s_{\underline u}(G)\}.
\end{equation}
One says that $H\leq G$ is {\em $\mathfrak{g}$-closed}  in $G$ if
$\mathfrak{g}_G(H)=H$, and $H$ is {\em $\mathfrak{g}$-dense} there
if $\mathfrak{g}_G(H)=G$.
Using this terminology, $\mathfrak g$ is the largest closure operator
with respect to which each {\em basic $\mathfrak{g}$-closed} set
$s_{\underline u}(G)$ is closed. (For details on closure operators
of topological groups, we refer the reader to the survey of 
Dikranjan \cite{DikCOTA}.)


For a category $\mathsf{T}$,  $\mathsf{Grp(T)}$ 
stands for the category of group objects in $\mathsf{T}$ and their group
homomorphisms, and $\mathsf{Ab(T)}$ for its full subcategory consisting
of the abelian group objects.
Thus, $\mathsf{Grp(Top)}$, $\mathsf{Grp(CompH)}$,
$\mathsf{Ab(Top)}$, and $\mathsf{Ab(CompH)}$ are the categories
of topological groups (without separation axioms), compact Hausdorff 
topological groups, and their abelian counterparts, respectively. 

For $\,G \in \mathsf{Grp(Top)}$, we denote by $\rho_G\colon G\rightarrow bG$
the {\em Bohr-compactification} of $G$, that is, 
the reflection of $G$ into $\mathsf{Grp(CompH)}$. One puts $G^+$ for
$G$ equipped with the initial topology induced by $\rho_G$.
The kernel of $\rho_G$ is denoted by $\mathbf{n}(G)$ and called the 
{\em von Neumann radical} of $G$. One says that 
$G$ is {\em maximally almost periodic} if $\mathbf{n}(G)=1$---in
other words, if $\rho_G$ is injective (cf.~\cite{NeuWig}).

\bigskip

In this paper, we study the closure operator $\mathfrak{g}$
(and its relatives that are introduced in section~\ref{sect:cont})
by means of modification of precompact group topologies.
Our main results in the context of the 
$\mathfrak{g}$-closure itself are Theorems 
\ref{thm:positive}-\ref{thm:gc:MAP} below, whose proofs appear in 
section~\ref{sect:gclosure}.

\begin{theorem} \label{thm:positive}
Let $H$ be a countable subgroup of a compact Hausdorff abelian 
topological group $K$. Then $H$ is $\mathfrak{g}$-closed in $K$.
\end{theorem}

Theorem~\ref{thm:positive} is a positive solution for two problems of 
Dikranjan,  Milan and Tonolo \cite[Problem~5.1, Question~5.2]{DikMilTon}, 
and it generalizes Theorem~\ref{fthm:BirDesSos}.
Independently, Dikranjan and Kunen \cite[1.9]{DikKun} 
and Mathias Beiglb\"ock, Christian Steineder and Reinhard Winkler
\cite{BeiSteWin} have also obtained Theorem~\ref{thm:positive}.

\begin{theorem} \label{thm:gc-realcompact}
Every $\mathfrak g$-closed subgroup of a compact Hausdorff
abelian group is realcompact.
\end{theorem}

\begin{theorem} \label{thm:gc:MAP}
For $G\in\mathsf{Ab(Top)}$, the following statements are equivalent:

\begin{list}{{\rm (\roman{enumi})}}
{\usecounter{enumi}\setlength{\labelwidth}{25pt}\setlength{\topsep}{-10pt}
\setlength{\itemsep}{-4pt} \setlength{\leftmargin}{25pt}}

\item
$G$ is maximally almost periodic;

\item
every countable subgroup of $G$ is $\mathfrak g$-closed;

\item
every cyclic subgroup of $G$ is $\mathfrak g$-closed.

\end{list}
\end{theorem}

Theorem~\ref{thm:gc:MAP} is a variant of the main result of
\cite{DikMilTon}, and so its novelty lies in its significantly simplified 
proof that uses techniques developed below.

\bigskip

One says that $G \in \mathsf{Grp(Top)}$ is {\em precompact} if for every 
neighborhood $U$ of the identity in $G$, there exists a finite subset 
$F\subseteq G$ such that $G=FU$. 

In order to describe  our method of investigating the 
closure operator $\mathfrak{g}$ (and its relatives that 
are introduced in section~\ref{sect:cont}), we first survey
some basic facts related to precompact abelian groups.
Fix $K \in \mathsf{Ab(CompH)}$, and 
let $A=\hat K$ be its Pontryagin dual. It is a 
well-known result of Comfort and Ross \cite{ComfRoss} that every
precompact Hausdorff group topology $\tau$ on $A$ coincides with
the initial topology on $A$ with respect to characters $\chi$ in a 
dense subgroup $H\leq K$. In other words, one has $\tau = \tau_H$,
where $\tau_H$ is the coarsest (group) topology on $A$ such that each
$\chi\colon A \rightarrow \mathbb{T}$ in $H$ is continuous. 
There is a close relationship between properties of $(A,\tau_H)$ and $H$.
To mention only two basic ones, the identity map  
$(A,\tau_{H_1})\rightarrow (A,\tau_{H_2})$ is continuous  (i.e., 
$\tau_{H_1} \supseteq \tau_{H_2}$) if and only if $H_1 \geq H_2$ 
(cf. \cite[2.6]{HernMaca}), and $(A,\tau_H)$ is metrizable if and only 
if $H$ is countable (cf. \cite[2.12]{HernMaca}). For a detailed
account of results of this nature, we refer the reader to \cite{HernMaca} and
\cite{ComfHernTrig}.

For $H\leq K$, consider
$H^\perp:= \bigcap\limits_{\chi\in H} \ker\chi$. It follows from the
Pontryagin duality that the dual of $A/H^\perp$ is the closure
$\bar H$ of $H$ in $K$, and $H^\perp = \bar H^\perp$. Thus,
although $H$ need not separate points in $A$, it always separates
the points of $A/H^\perp$. In particular, $H$ separates the points 
of $A$ if and only if $H$ is dense 
in $K$ (cf. \cite{ComfRoss}).
Therefore, in the context of characters
and duality, the precompact (but not necessarily Hausdorff) group 
$(A,\tau_H)$ is equivalent to the precompact Hausdorff group
$(A/H^\perp,\tau_H)$, because the latter is the maximal Hausdorff
quotient of the earlier (see \cite[2.10]{HernMaca}).

The key to our approach is our viewing $\mathfrak{g}$ as a functor
that assigns to each precompact abelian group $(A,\tau_H)$
the group $(A,\tau_{\mathfrak g_{\empty_K}(H)})$. We show that
$\mathfrak g_K(H)$ consists of the elements of $K$ that are 
sequentially continuous as characters of $(A,\tau_H)$
(i.e., $\chi(a_n)\longrightarrow \chi(a_0)$ 
for every sequence $\{a_n\}$ in $A$ such that
$a_n \longrightarrow a_0$ in $\tau_H$; cf. Lemma~\ref{lemma:f}(b)).

\bigskip

In the course of investigating $\mathfrak{g}$, we discovered
that many of its properties remain true if sequences are replaced
with  a ``nice" method of choosing filter-bases in $\hat G$. 
Thus, in order to achieve a greater generality, in
section~\ref{sect:cont}, we study such mutants of 
$\mathfrak{g}$. Each closure operator of this type gives rise to 
a full coreflective subcategory of the category of precompact abelian groups.
In section~\ref{sect:realcompact}, we present such a ``mutant" 
$\mathfrak{l}$ of $\mathfrak{g}$ (in the sense of 
section~\ref{sect:cont}) that coincides with the 
$G_\delta$-closure on compact Hausdorff abelian groups
(cf.~Theorem~\ref{thm:char:l}); thus, 
$\mathfrak{l}$-closed subgroups are the realcompact ones, 
and $\mathfrak{l}$-dense subgroups are the dense pseudocompact ones.
Finally, in section~\ref{sect:gclosure},  the general theory developed in 
section~\ref{sect:cont} is applied to the closure operator $\mathfrak{g}$,
and Theorems~\ref{thm:positive}-\ref{thm:gc:MAP} are proven.

\section{Continuity types}

\label{sect:cont}

In this section, we develop the filter-base analogue of the 
$\mathfrak{g}$-closure. Let $G \in \mathsf{Ab(Top)}$ and let $\mathcal{F}$ 
be a filter-base on its Pontryagin dual $\hat G$. We put
\begin{equation}
z_{\mathcal{F}}(G) = \{ \chi\in G \mid \chi(\mathcal{F}) \longrightarrow 0 
\text{ in } \mathbb{T} \},
\end{equation}
where elements of $G$ are seen as characters of $\hat G$, and thus 
they are denoted by Greek letters.

\begin{lemma} \label{lemma:fCLOP:funct}
Let $\varphi\colon G_1 \rightarrow G_2$ be a morphism in 
$\mathsf{Ab(Top)}$, and
$\mathcal{F}$ be a filter-base on $\hat G_2$. Then
\begin{equation}
z_{\hat \varphi (\mathcal{F})}(G_1) = \varphi^{-1}(z_{\mathcal{F}}(G_2)).
\end{equation}
\end{lemma}

\begin{proof}
Since $\hat \varphi\colon \hat G_2 \rightarrow \hat G_1$ is defined by
$\hat \varphi(x)=x \circ \varphi$, 
\begin{align}
\chi\in z_{\hat \varphi (\mathcal{F})}(G_1) & \Longleftrightarrow 
\chi(\hat\varphi (\mathcal{F})) \longrightarrow 0 \Longleftrightarrow 
\varphi(\chi)(\mathcal{F})\longrightarrow 0\Longleftrightarrow 
\varphi(\chi) \in z_{\mathcal{F}}(G_2),
\end{align}
as desired.
\end{proof}

A {\em continuity type} $\mathcal{T}$ on a subcategory $\mathsf{A}$
of $\mathsf{Grp(Top)}$
is an assignment of a collection $\mathcal{T}(G)$ of filter-bases to 
each $G \in \mathsf{A}$ in a way that 
$\varphi(\mathcal{F})\in \mathcal{T} (G_2)$ for 
every morphism $\varphi\colon G_1 \rightarrow G_2$ in $\mathsf{A}$
and every $\mathcal{F}\in\mathcal{T}(G_1)$.  
Examples~\ref{example:S} and~\ref{example:C} 
below explain the motivation behind the term ``continuity type."

\begin{proposition} \label{prop:f-closure}
Let $\mathcal{T}$ be a continuity type on $\mathsf{Ab(Top)}$.
For every $G \in \mathsf{Ab(Top)}$ and $H \leq G$, put
\begin{equation}
\mathfrak{f}_G(H)=\bigcap\{z_\mathcal{F}(G)\mid\mathcal{F}\in\mathcal{T}(\hat G),
H \leq z_\mathcal{F}(G)\}.
\end{equation}
Then $\mathfrak{f}$ is an idempotent
closure operator on $\mathsf{Ab(Top)}$.
\end{proposition}

We call $\mathfrak{f}$ the closure operator 
{\em associated with $\mathcal{T}$}.

\begin{proof}
Clearly, $\mathfrak{f}$ is extensive and monotone, so we prove that
each continuous homomorphism $\varphi\colon G_1 \rightarrow G_2$ 
is $\mathfrak{f}$-continuous. Let $H \leq G_2$.
Since $\hat \varphi \colon \hat G_2 \rightarrow \hat G_1$ is continuous,
$\hat\varphi(\mathcal{F}_2)\in \mathcal{T}(\hat G_1)$ for every
$\mathcal{F}_2\in \mathcal{T}(\hat G_2)$. Thus, 
\begin{align}
\mathfrak f_{G_1}(\varphi^{-1}(H)) & = 
\bigcap\{z_{\mathcal{F}_1}(G_1)\mid\mathcal{F}_1\in\mathcal{T}(\hat G_1),
\varphi^{-1}(H) \leq z_{\mathcal{F}_1}(G_1)\} \\
&  \subseteq \bigcap\{z_{\hat\varphi(\mathcal{F}_2)}(G_1)\mid
\hat\varphi(\mathcal{F}_2) \in\mathcal{T}(\hat G_1),
\varphi^{-1}(H) \leq z_{\hat\varphi (\mathcal{F}_2)}(G_1)\} \\
& \stackrel  {\ref{lemma:fCLOP:funct}} =
\bigcap\{\varphi^{-1}(z_{\mathcal{F}_2}(G_2))\mid
\mathcal{F}_2 \in\mathcal{T}(\hat G_2),
(H) \leq z_{\mathcal{F}_2}(G_2)\} = 
\varphi^{-1} (\mathfrak f_{G_2}(H)),
\end{align}
as desired.  Obviously, $\mathfrak{f}$ is idempotent 
(cf. \cite[Exercise~2.D(a)]{DT}).
\end{proof}

A continuity type $\mathcal{T}$ is {\em regular on $G$} if 
$\mathcal{T}(G)=\mathcal{T}(G_d)$, where $G_d$ is the group $G$ with
discrete topology; a subgroup $N$ of $G \in \mathsf{Ab(Top)}$ is 
{\em $\mathcal{T}$-dually embedded} into $G$ if for every
$\mathcal{F}_2 \in \mathcal{T}(\hat N)$, there is 
$\mathcal{F}_1 \in \mathcal{T}(\hat G)$ such that
$\hat i_N(\mathcal{F}_1)=\mathcal{F}_2$.

If $G \in \mathsf{Ab(Top)}$, then $G^+$ is the group $G$ with the 
topology of pointwise convergence on $\hat G$, and
$\mathbf{n}(G)=\bigcap\limits_{x \in \hat G} \ker x$.
It follows from the universal
property of $bG$ that $\widehat{bG}=\hat G_d$.

\begin{proposition} \label{prop:contype}
Let $\mathcal{T}$ be a continuity type on $\mathsf{Ab(Top)}$, 
$\mathfrak{f}$ its associated closure operator, and
$G \in \mathsf{Ab(Top)}$.

\begin{list}{{\rm (\alph{enumi})}}
{\usecounter{enumi}\setlength{\labelwidth}{25pt}\setlength{\topsep}{-10pt}
\setlength{\itemsep}{-4pt} \setlength{\leftmargin}{25pt}}

\item
One has $\mathbf{n}(G)\subseteq \mathfrak f_G(\{0\})$.

\item
If every principal ultrafilter on $\hat G$ is in $\mathcal{T}(\hat G)$,
then $\mathbf{n}(G)=\mathfrak f_G(\{0\})$.

\item
If \,$\mathcal{T}$ is regular on $\hat G$ , then 
$\mathfrak f_G(H) = \mathfrak f_{G^+}(H)$ for every $H\leq G$.

\item
If\, $\mathcal{T}$ is regular, then 
every dense subgroup of\, $G$ is $\mathcal{T}$-dually embedded
into $G$.

\item
If $N$ is $\mathcal{T}$-dually embedded into
$G$, then $\mathfrak{f}_N(H)=N\cap\mathfrak{f}_G(H)$ for
every $H\leq N$.

\end{list}
\end{proposition}

\begin{proof}
(a) If $\chi \in \mathbf{n}(G)$, then $\chi(\mathcal{F}) = \{ \{0\} \}$ for
every filter-base $\mathcal{F}$ on $\hat G$, and so $\chi\in z_\mathcal{F}(G)$. 
Thus, $\mathbf{n}(G)\subseteq z_\mathcal{F}(G)$ for every
$\mathcal{F}$, and therefore 
$\mathbf{n}(G)\subseteq \bigcap\limits_{\mathcal{F}\in\mathcal{T}(\hat G)}
z_\mathcal{F}(G)=\mathfrak f_G(\{0\})$.

(b) Let $x \in \hat G$. By the assumption, the principal ultrafilter
$\dot x \in \mathcal{T}(\hat G)$, and $z_{\dot {x}}(G) = \ker x$. 
Thus, $\mathfrak f_G(\{0\}) \subseteq \ker x$. This holds for 
every $x \in \hat G$, and therefore (using (a))
$\mathbf{n}(G)=\mathfrak f_G(\{0\})$.

(c) The continuous homomorphisms $G \rightarrow \hat G \rightarrow bG$ give 
rise to  $\hat G_d=\widehat{bG}\rightarrow\widehat{G^+}\rightarrow\hat G$, and 
thus $\mathcal{T}(\hat G_d)\subseteq\nolinebreak\mathcal{T}(\widehat{G^+})
\subseteq\mathcal{T}(\hat G)$. Since $\mathcal{T}$ is regular on 
$\hat G$, the two ends of this inclusion are equal, and therefore
$\mathcal{T}(\widehat{G^+})=\mathcal{T}(\hat G)$.
The definition of $z_\mathcal{F}$ is independent of the topology,
and thus $z_\mathcal{F}(G)=\nolinebreak z_\mathcal{F}(G^+)$ for every 
filter-base $\mathcal{F}$ on the \underline{set} $\hat G=\widehat{G^+}$. 
Hence, the collections whose intersection give 
$\mathfrak f_G(H)$ and $\mathfrak f_{G^+}(H)$ coincide.

(d) Let $D \leq G$ be a dense subgroup. Then
$\hat G$ and $\hat D$ have the same underlying group, and so
$\mathcal{T}(\hat D)=\mathcal{T}(\hat D_d)=\mathcal{T}(\hat G_d)=
\mathcal{T}(\hat G)$, because $\mathcal{T}$ is regular.

(e) Clearly, $\mathfrak{f}_N(H)\subseteq \mathfrak{f}_G(H)$, because
the inclusion $N \rightarrow G$ is $\mathfrak f$-continuous. Let
$\mathcal{F}_2 \in \mathcal{T}(\hat N)$, and let 
$\mathcal{F}_1\in\mathcal{T}(\hat G)$ be such that 
$\hat i_N(\mathcal{F}_1)=\mathcal{F}_2$ (the existence of
$\mathcal{F}_1$ is guaranteed as $N$ is $\mathcal{T}$-dually embedded
into $G$). By Lemma~\ref{lemma:fCLOP:funct},  $z_{\mathcal{F}_2}(N)=
z_{\hat i_{\empty_N} (\mathcal{F}_1)}(N)=i_N^{-1}(z_{\mathcal{F}_1}(G))=
N \cap z_{\mathcal{F}_1}(G)$, so if $H \leq z_{\mathcal{F}_2}(N)$,
then $z_{\mathcal{F}_2}(N)\supseteq N \cap \mathfrak{f}_G(H)$.
Therefore, $\mathfrak{f}_N(H)\supseteq  N\cap \mathfrak{f}_G(H)$,
as desired.
\end{proof}

\begin{theorem} \label{thm:reduce:bohr}
Let $\mathcal{T}$ be a regular continuity type on $\mathsf{Ab(Top)}$ and
$\mathfrak{f}$ its associated closure operator. Then
$\mathfrak{f}_G(H) = \rho^{-1}_G(\mathfrak{f}_{bG}(\rho_G(H)))$
for every abelian topological group $G$ and subgroup $H \leq G$.
\end{theorem}

\begin{proof}
Set $b^+G=\rho_G(G)$; it is a dense subgroup of $bG$. Observe that
$bG = b(G^+)$ (and so $b^+(G^+)=b^+ G$), and $b^+G^+=G^+/\mathbf{n}(G)$.
Since $b^+(G^+)$ is dense in $b(G^+)$,
by Proposition~\ref{prop:contype}(d), $b^+(G^+)$ is 
$\mathcal{T}$-dually embedded into $bG$. Thus, one has
\begin{align}
\mathfrak{f}_G(H)
\stackrel{\rm \ref{prop:contype}(c)} = 
\mathfrak{f}_{G^+}(H)
& \hspace{6pt} =\hspace{6pt}
\rho^{-1}_{G^+}(\mathfrak{f}_{b^+(G^+)}(\rho_{G^+}(H)))  
\label{eq:exp:1}\\
& \stackrel {\rm \ref{prop:contype}(e)}=
\rho^{-1}_G(\mathfrak{f}_{bG}(\rho_G(H))\cap b^+G) 
 = \rho^{-1}_G(\mathfrak{f}_{bG}(\rho_G(H))).
\end{align}
Note that the last equality in (\ref{eq:exp:1}) holds because
every continuous character of $G^+$ factors through 
$\rho_{G^+}$ (so $\widehat{G^+}=\widehat{b^+(G^+)}$ as sets), 
and $\mathcal{T}$ is regular.
\end{proof}

\begin{example}
For every topological group $G$, put $\mathcal{D}(G)$ to be the set of all
filter-bases on $G$. The associated closure operator is 
discrete on every $K \in \mathsf{Ab(CompH)}$, because for every 
(not necessary closed)
subgroup $H$, there is a filter $\mathcal{F}$ on $A=\hat K$
such that $H=z_\mathcal{F}(K)$ (cf.~\cite[2.1]{BeiSteWin}). 
Therefore, by Theorem~\ref{thm:reduce:bohr},
for $G \in \mathsf{Ab(Top)}$,
$\mathfrak{f}^\mathcal{D}_G(H)=\rho^{-1}_G(\rho_G(H))=H+\mathbf{n}(G)$. Hence,
$\mathfrak{f}^\mathcal{D}_G(H) = c_\mathbf{n}$, the minimal closure operator
associated with $\mathbf{n}$ (cf.~\cite[3.5-3.6]{DikCOTA}).
\end{example}

Let $\mathsf{A}$ be a subcategory of $\mathsf{Grp(Top)}$ and $\mathcal{T}$ be
a continuity type on $\mathsf{A}$. For $G_1,G_2 \in \mathsf{A}$,
we say that a (not necessarily continuous)
homomorphism $\varphi\colon G_1 \rightarrow G_2$ is 
{\em $\mathcal{T}$-continuous} if $\varphi(\mathcal{F})\longrightarrow e$ 
in $G_2$ for every $\mathcal{F}\in\mathcal{T}(G_1)$ 
such that $\mathcal{F}\longrightarrow e$ in $G_1$. 

\begin{example} \label{example:S}
For every topological group $G$, put $\mathcal{S}(G)$ to be the set of all
filter-bases on $G$ generated by sequences (i.e., of the form 
$\{\{a_n \mid n \geq m\}\mid m \in\mathbb{N}\}$, where $\{a_n\}\subseteq G$).
A homomorphism $\varphi\colon G_1 \rightarrow G_2$ is $\mathcal{S}$-continuous
if and only if it is  sequentially continuous.
For $G\in\mathsf{Ab(Top)}$ and $\underline u \in \hat G^\mathbb{N}$,
one has $s_{\underline u}(G)=z_\mathcal{F}(G)$, where 
$\mathcal{F}\in\mathcal{S}(G)$ is the filter base generated by the 
sequence $\underline u$, and vice versa, to each 
$\mathcal{F}\in\mathcal{S}(G)$ there is a corresponding sequence that generates
it. Therefore, the associated closure operator of $\mathcal{S}$ is 
precisely $\mathfrak{g}$, which was defined earlier.
\end{example}

A map $f\colon X \rightarrow Y$ between topological spaces is said to be
{\em countably continuous} if $f$ is continuous on every
countable subset of $X$.

\begin{example} \label{example:C}
For every topological group $G$, put $\mathcal{C}(G)$ to be the set
of filter-bases on $G$ that contain a countable subset of $G$. 
In Lemma~\ref{lemma:l:cc}, we show that a homomorphism 
$\varphi\colon G_1 \rightarrow G_2$ is $\mathcal{C}$-continuous 
if and only if it is countably continuous. Moreover, in 
Theorem~\ref{thm:char:l}, we 
prove that the closure operator $\mathfrak{l}$ associated with 
$\mathcal{C}$  coincides with the $G_\delta$-closure on precompact 
Hausdorff abelian groups.
\end{example}

\begin{convention}
Until the end of the section, $\mathcal{T}$ is 
a fixed regular continuity type on $\mathsf{Ab(Top)}$,
and $\mathfrak{f}$ is its associated closure operator.
\end{convention}

Motived by Theorem~\ref{thm:reduce:bohr},  we confine 
our attention to $\mathsf{Ab(CompH)}$.
In the sequel, we put $\mathcal{F} \stackrel H \longrightarrow a_0$
whenever a filter-base $\mathcal{F}$
converges to a point $a_0$ in $A$ in the topology $\tau_H$.

\begin{samepage}
\begin{lemma} \label{lemma:f}
Let $K \in \mathsf{Ab(CompH)}$, $A=\hat K$, and $H \leq K$. Put
$\mathcal{T}_H(A) = \{\mathcal{F}\in\mathcal{T}(A)\mid 
\mathcal{F} \stackrel H \longrightarrow 0 \}$. Then:

\begin{list}{{\rm (\alph{enumi})}}
{\usecounter{enumi}\setlength{\labelwidth}{25pt}\setlength{\topsep}{-10pt}
\setlength{\itemsep}{-4pt} \setlength{\leftmargin}{25pt}}

\item
$\mathfrak f_K(H)=\{ \chi \in K\mid \chi(\mathcal{F}) \longrightarrow 0
\text{ for all } \mathcal{F}\in \mathcal{T}_H(A) \}$;

\item
$\mathfrak f_K(H)=\{ \chi \in K\mid \chi\colon(A,\tau_H) 
\rightarrow \mathbb{T} \text{ is $\mathcal{T}$-continuous}\}$;

\item
the identity homomorphism $(A,\tau_H)\rightarrow 
(A,\tau_{\mathfrak{f}_{\empty_K}(H)})$ is $\mathcal{T}$-continuous.

\end{list}
\end{lemma}
\end{samepage}

\begin{proof}
(a) For $\mathcal{F}\in \mathcal{T}(A)$, $H \leq z_\mathcal{F}(K)$ if and only if
$\chi(\mathcal{F}) \longrightarrow 0$ for every $\chi \in H$---in other words,
$\mathcal{F} \stackrel H \longrightarrow 0$, which is the same as
$\mathcal{F}\in\mathcal{T}_H(A)$. Therefore,
\begin{align}
\mathfrak{f}_K(H)=\bigcap\{ z_\mathcal{F}(K) \mid 
\mathcal{F}\in \mathcal{T}_H(A)\} =
\{\chi \in K \mid \chi(\mathcal{F})\longrightarrow 0 \text{ for all }
\mathcal{F} \in \mathcal{T}_H(A)\}.
\end{align}

(b) Let $\chi \in \mathfrak{f}_K(H)$, and suppose that 
$\mathcal{F} \in \mathcal{T}(A)$ is such that $\mathcal{F}\longrightarrow 0$
in $\tau_H$. Then $\mathcal{F}\in\nolinebreak\mathcal{T}_H(A)$, and so
$\chi(\mathcal{F})\longrightarrow 0$ by (a). Thus,
$\chi$ is $\mathcal{T}$-continuous. Conversely, 
suppose that $\chi\colon (A,\tau_H)\rightarrow \mathbb{T}$ is 
$\mathcal{T}$-continuous, and let $\mathcal{F}\in \mathcal{T}_H(A)$.
Then $\chi(\mathcal{F})\longrightarrow 0$, because $\chi$ is 
$\mathcal{T}$-continuous, and therefore $\chi\in \mathfrak{f}_K(H)$
by (a).

(c) If $\mathcal{F}\in \mathcal{T}_H(A)$, then by (a),
$\chi(\mathcal{F}) \longrightarrow 0$ for every 
$\chi \in \mathfrak{f}_K(H)$, and therefore
$\mathcal{F} \stackrel {\mathfrak{f}(H)}\longrightarrow 0$.
\end{proof}

The next theorem is an immediate consequence of Lemma~\ref{lemma:f}(b).

\begin{theorem} \label{thm:f}
Let $K \in \mathsf{Ab(CompH)}$ and $A=\hat K$.
A subgroup $H \leq K$ is $\mathfrak{f}$-closed in $K$ if and only if every 
$\mathcal{T}$-continuous character of $(A,\tau_H)$ is continuous.
\qed
\end{theorem}

\begin{proposition} \label{prop:f-funct}
Let $A$ and $B$ be discrete abelian groups, $K=\hat A$ 
and $C=\hat B$ their Pontryagin duals, $H \leq K$, and $L\leq C$. For every
homomorphism $\varphi\colon A \rightarrow B$,
the following statements are equivalent:

\begin{list}{{\rm (\roman{enumi})}}
{\usecounter{enumi}\setlength{\labelwidth}{25pt}\setlength{\topsep}{-10pt}
\setlength{\itemsep}{-4pt} \setlength{\leftmargin}{25pt}}

\item $\varphi\colon (A,\tau_H) \rightarrow (B,\tau_L)$ is 
$\mathcal{T}$-continuous;

\item $\chi\circ\varphi$ is $\mathcal{T}$-continuous on $(A,\tau_H)$ for
every $\mathcal{T}$-continuous character $\chi$ of $(B,\tau_L)$;

\item
$\mathfrak{f}\varphi\colon(A,\tau_{\mathfrak{f}_{\empty_K}(H)})\rightarrow 
(B,\tau_{\mathfrak{f}_{\empty_C}(L)})$ is continuous.
\end{list}
\end{proposition}

\begin{proof}
Clearly, (i) $\Rightarrow$ (ii), because the composite of two
$\mathcal{T}$-continuous homomorphisms is $\mathcal{T}$-continuous. 
The implication (ii) $\Rightarrow$ 
(iii) follows from Lemma~\ref{lemma:f}(b), because (ii) states that
$\chi\circ\varphi \in \mathfrak{f}_K(H)$ for every
$\chi \in \mathfrak{f}_C(L)$. To show (iii) $\Rightarrow$ (i), notice that
$\varphi$ is the composite of\linebreak 
$(A,\tau_H)\rightarrow (A,\tau_{\mathfrak{f}_{\empty_K}(H)})$, 
$\mathfrak f \varphi$, 
and $(B,\tau_{\mathfrak{f}_{\empty_C}(L)}) \rightarrow (B,\tau_L)$. Since these  
three maps are all $\mathcal{T}$-continuous (cf. Lemma~\ref{lemma:f}(c)), 
their composite is also $\mathcal{T}$-continuous.
\end{proof}

Inspired by the notions of $k$-group and $s$-group introduced by Noble 
\cite{Noble3} \& \cite{Noble}, and the term ``$kk$-group" used by Deaconu 
\cite{PHDDeaconu},
we say that a topological group $G$ is a {\em $\mathcal{T}k$-group} 
if every $\mathcal{T}$-continuous homomorphism of $G$ into a 
compact Hausdorff group is continuous.

\pagebreak[2]

\begin{proposition} \label{prop:Tk-grp}
Let $K \in \mathsf{Ab(CompH)}$, $A=\hat K$, and $H \leq K$.
The following are equivalent:

\begin{list}{{\rm (\roman{enumi})}}
{\usecounter{enumi}\setlength{\labelwidth}{25pt}\setlength{\topsep}{-10pt}
\setlength{\itemsep}{-4pt} \setlength{\leftmargin}{25pt}}

\item
every $\mathcal{T}$-continuous homomorphism 
$\varphi\colon (A,\tau_H) \rightarrow (B,\tau_L)$ into a precompact 
abelian group (where $B$ is a discrete group and $L \leq \hat B$) 
is continuous;

\item
$(A,\tau_H)$ is an $\mathcal{T}k$-group;

\item
every $\mathcal{T}$-continuous character
$\chi\colon (A,\tau_H) \rightarrow \mathbb{T}$ is continuous.

\end{list}
\end{proposition}

\begin{proof}
(i) $\Rightarrow$ (ii) $\Rightarrow$ (iii) is obvious, because
every compact group is precompact, and $\mathbb{T}$ is compact.
In order to show (iii) $\Rightarrow$ (i), let
$\varphi\colon (A,\tau_H) \rightarrow (B,\tau_L)$  be an 
$\mathcal{T}$-continuous
homomorphism. By Proposition~\ref{prop:f-funct}, 
$\chi \circ \varphi$ is $\mathcal{T}$-continuous on $(A,\tau_H)$ for 
every character $\chi \in L$ of $(B,\tau_L)$. Applying (iii) to 
the character $\chi \circ \varphi$ yields that it is continuous, in other
words, $\chi \circ \varphi \in H$. Therefore, $\chi\circ\varphi \in H$ 
for every $\chi \in L$, and hence $\varphi$ is continuous.
\end{proof}

A combination of Theorem~\ref{thm:f} and Proposition~\ref{prop:Tk-grp} 
yields:

\begin{theorem} \label{thm:Tk-closed}
Let $K \in \mathsf{Ab(CompH)}$ and $A=\hat K$.
A subgroup $H\leq K$ is $\mathfrak{f}$-closed if and only if\,
$(A,\tau_H)$ is an $\mathcal{T}k$-group. In particular,
$(A,\tau_{\mathfrak f_{\empty_K}(H)})$ is an $\mathcal{T}k$-group for every
$H \leq K$.
\qed
\end{theorem}

\bigskip

Theorem~\ref{thm:subcat} below is modeled on the relationship between the 
category of Hausdorff topological spaces with $k$-continuous maps,
Hausdorff topological spaces with continuous maps, and Hausdorff $k$-spaces
(cf.~\cite[7.2]{Bor2}).
In what follows, $\mathsf{Pr_\mathcal{T}Ab}$ is the category of precompact abelian 
groups and their $\mathcal{T}$-continuous homomorphisms, $\mathsf{PrAb}$ stands 
for the category of precompact abelian groups and their continuous homomorphisms,
and $\mathsf{\mathcal{T}kPrAb}$ is the full
subcategory of $\mathsf{PrAb}$ consisting of the $\mathcal{T}k$-groups.

\begin{theorem} \label{thm:subcat}
The assignment
\begin{align*}
\mathfrak f \colon \mathsf{Pr_\mathcal{T}Ab} & 
\longrightarrow \mathsf{\mathcal{T}kPrAb} \\
\hspace{120pt}
(A,\tau_H) & \longmapsto (A,\tau_{\mathfrak f_{\empty_K} (H)}), 
\quad (K=\hat A,H\leq K)
\end{align*}
is a functor. Furthermore,

\begin{list}{{\rm (\alph{enumi})}}
{\usecounter{enumi}\setlength{\labelwidth}{25pt}\setlength{\topsep}{-10pt}
\setlength{\itemsep}{-4pt} \setlength{\leftmargin}{25pt}}

\item 
$\mathfrak f \colon \mathsf{Pr_\mathcal{T}Ab}  \longrightarrow 
\mathsf{\mathcal{T}kPrAb}$ is  an equivalence of categories;

\item
$\mathfrak f \colon \mathsf{PrAb}  \longrightarrow 
\mathsf{\mathcal{T}kPrAb}$ is a coreflection.
\end{list}
\end{theorem}

\begin{proof}
It was explained earlier that every precompact abelian group has the form
$(A,\tau_H)$, where $A$ is a discrete abelian group and  $H \leq \hat A$
($A$ and $H$ are uniquely determined). By Theorem~\ref{thm:Tk-closed},
$\mathfrak{f} (A,\tau_H) = (A,\tau_{\mathfrak{f}_{\empty_K}(H)}) 
\in \mathsf{\mathcal{T}kPrAb}$,
and by Proposition~\ref{prop:f-funct}, $\mathfrak{f}\varphi$ is continuous 
for every $\mathcal{T}$-continuous homomorphism 
$\varphi\colon (A,\tau_H) \rightarrow (B,\tau_L)$.
Therefore, $\mathfrak{f}$ is well-defined, and it is a functor.

(a) Clearly, $\mathfrak{f}$ and the inclusion 
$I\colon \mathsf{\mathcal{T}kPrAb} \longrightarrow \mathsf{Pr_\mathcal{T}Ab}$ are
faithful, and by Proposition~\ref{prop:f-funct}, they are also full. Since 
$\mathfrak{f}$ is idempotent, 
$\mathfrak{f}I=\operatorname{Id}_{\mathsf{\mathcal{T}kPrAb}}$.
By Lemma~\ref{lemma:f}(c), $(A,\tau_H)$ and 
$(A,\tau_{\mathfrak{f}_{\empty_K}(H)})$
are isomorphic in $\mathsf{Pr_\mathcal{T}Ab}$, and therefore
$I \mathfrak{f} \cong \operatorname{Id}_{\mathsf{Pr_\mathcal{T}Ab}}$. Hence,
$\mathfrak{f}$ is an equivalence of categories.

(b) For $(B,\tau_L)\in \mathsf{PrAb}$ (where $B$ is a discrete abelian group
and $L \leq C:=\hat B$),
the counit $\gamma_{(B,\tau_L)}\colon (B,\tau_{\mathfrak f_{\empty_C}(L)}) 
\rightarrow (B,\tau_L)$ is continuous (because $L \leq \mathfrak f_C(L)$).
By Proposition~\ref{prop:f-funct}, every continuous
homomorphism $\varphi\colon (A,\tau_H) \rightarrow (B,\tau_L)$ 
gives rise to a continuous homomorphism 
$\mathfrak f \varphi \colon (A,\tau_{\mathfrak{f}_{\empty_K}(H)}) 
\rightarrow (B,\tau_{\mathfrak f_{\empty_C}(L)})$. Whenever $(A,\tau_H)$ is 
an $\mathcal{T}k$-group, $\mathfrak{f}_K(H)=H$ (cf. Theorem~\ref{thm:Tk-closed}), 
and therefore $\varphi$ factors through 
$\gamma_{(B,\tau_L)}$ uniquely.
\end{proof}


\section{Application I: The $\mathbf{G_\delta}$-closure}

\label{sect:realcompact}

Recall that for a topological space $X$, the {\em $G_\delta$-closure} 
$\bar A^\delta$ of $A \subseteq X$ consists of all  points $x \in X$ such 
that each $G_\delta$-set that contains $x$ meets $A$; $A$ is {\em
$G_\delta$-dense} in $X$ if $\bar A ^\delta =X$.

Let $\mathcal{C}$ be the continuity type described in 
Example~\ref{example:C}, and set $\mathfrak{l}$ to be 
its associated closure operator. In this section, we study the 
closure operator $\mathfrak l$ on precompact Hausdorff abelian groups.

\begin{theorem} \label{thm:char:l}
Let $P$ be a precompact Hausdorff abelian group. Then
$\mathfrak l_P(H) = \bar H^\delta$ for every subgroup $H \leq P$. In other 
words, $\mathfrak l$ coincides with the $G_\delta$-closure on precompact 
Hausdorff abelian groups.
\end{theorem}

Before proving Theorem~\ref{thm:char:l}, we present some basic properties 
of $\mathfrak{l}$.
Recall that a subgroup $N$ of $G\in\mathsf{Ab(Top)}$ is said to be
{\em dually embedded} if $\hat i_N\colon \hat G \rightarrow \hat N$
is surjective, where $i_N$ is the inclusion.

\begin{proposition} \label{prop:l:basic} 
For every $G\in\mathsf{Ab(Top)}$:

\begin{list}{{\rm (\alph{enumi})}}
{\usecounter{enumi}\setlength{\labelwidth}{25pt}\setlength{\topsep}{-10pt}
\setlength{\itemsep}{-4pt} \setlength{\leftmargin}{25pt}}

\item
$\mathfrak{l}_G(\{0\})=\mathbf{n}(G)$;

\item
$\mathfrak{l}_G(H)=\mathfrak{l}_{G^+}(H)$ for every $H\leq G$;

\item
if $N\leq G$ is dually embedded in $G$, then 
$\mathfrak{l}_N(H)=\mathfrak{l}_G(H)\cap N$ for every $H \leq N$.

\end{list}
\end{proposition}

\begin{proof}
(a) and (b) follow from Proposition~\ref{prop:contype} (b) and (c), respectively.
By Proposition~\ref{prop:contype}(e), in order to show (c), it
suffices to show that $N$ is $\mathcal{C}$-dually embedded in $G$.
To that end, let $\mathcal{F}_2 \in \mathcal{C}(\hat N)$,
and let $D_2 \in \mathcal{F}_2$ be a countable subset.
Since $\hat i_N\colon \hat G \rightarrow \hat N$ is surjective, 
there exists a countable subset $D_1$ of $\hat G$ such that 
$\hat i(D_1)=D_2$. It is easy to see that for 
\begin{align}
\mathcal{F}_1=\{\hat i_N^{-1}(F) \mid F \in \mathcal{F}_2\} \cup
\{\hat i_N^{-1}(F)\cap D_1 \mid F \in \mathcal{F}_2\},
\end{align}
$\mathcal{F}_1 \in \mathcal{C}(\hat G)$, and 
$\hat i_N(\mathcal{F}_1)=\mathcal{F}_2$.
\end{proof}

It follows by Pontryagin duality that
every subgroup of a compact Hausdorff abelian group
is dually embedded, and thus Proposition~\ref{prop:l:basic}(c)
yields:

\begin{corollary} \label{cor:l:her}
For every $K\in \mathsf{Ab(CompH)}$ and subgroups
$H\leq N \leq K$,
$\mathfrak{l}_N(H)=\mathfrak{l}_K(H)\cap N$ holds.
In other words, $\mathfrak{l}$ is hereditary on
$\mathsf{Ab(CompH)}$, and thus $\mathfrak{l}$ is hereditary on every 
precompact Hausdorff abelian group.
\qed
\end{corollary}

\begin{lemma} \label{lemma:l:cc}
A homomorphism $\varphi\colon G_1 \rightarrow G_2$ 
is $\mathcal{C}$-continuous if and only if it is 
countably continuous.
\end{lemma}

\begin{proof}
Let $D$ be a countable subset of $G_1$, and set $S=\langle D \rangle$,
the (countable) subgroup generated by $D$. Every filter-base
$\mathcal{F}$ on $S$ belongs to $\mathcal{C}(G_1)$. Thus,  if
$\mathcal{F}\longrightarrow e$ in $S$, then
$\varphi(\mathcal{F}) \longrightarrow e$ in $G_2$, because
$\varphi$ is $\mathcal{T}$-continuous. Therefore, $\varphi\rest_S$ is 
continuous, and in particular $\varphi\rest_D$ is continuous.

Conversely, suppose that $\varphi$ is countably continuous, and
let $\mathcal{F}\in \mathcal{C}(G_1)$ such that 
$\mathcal{F} \longrightarrow e$ in $G_1$. Let $D\in \mathcal{F}$ be countable.
By replacing $\mathcal{F}$ with the filter-base 
$\{F \in \mathcal{F} \mid F \subseteq D\}$, we may assume that every set in
$\mathcal{F}$ is contained in $D$. Set $D_0=D\cup\{e\}$. Since 
$\varphi$ is countably continuous,  $\varphi\rest_{D_0}$ is continuous. 
Thus, $\varphi\rest_{D_0}(\mathcal{F})\longrightarrow e$, because
$\mathcal{F}\longrightarrow e$ in $D_0$. 
\end{proof}

\begin{proof}[Proof of Theorem~\ref{thm:char:l}.] 
Since both $\mathfrak l$ and the $G_\delta$-closure are hereditary on 
precompact Hausdorff groups (cf. Corollary~\ref{cor:l:her}), it suffices 
to show the theorem for the case where $P=K$ is compact. 
Set $A=\hat K$. By Lemma~\ref{lemma:f}(b), $\chi \in \mathfrak l_K(H)$ if 
and only if $\chi$ is $\mathcal C$-continuous with respect to $\tau_H$. By 
Lemma~\ref{lemma:l:cc}, this amounts to $\chi$ being countably continuous 
on $(A,\tau_H)$. We follow the idea of \cite[3.1]{ComfHernTrig} in order 
to show that the last statement is equivalent to $\chi \in \bar H^\delta$.

Since $K \in \mathsf{Ab(CompH)}$, $K$ embeds into the product
$\mathbb{T}^A$. For each $(t_a)_{a\in A} \in \mathbb{T}^A$ and countable 
subset $C \subseteq A$, the set
$\{(s_a)_{a\in A}\mid s_c=t_c \mbox{ for all } c \in C \}$
is $G_\delta$ in $\mathbb{T}^A$ (because $\mathbb{T}$ has a countable 
pseudocharacter), and these sets form a base for the $G_\delta$-topology 
on $\mathbb{T}^A$. Thus, for every countable subset $C \subseteq A$,
\begin{equation}
N(\chi,C)=\{\psi \in K \mid \psi\rest_C = \chi\rest_C\}
\end{equation}
is a basic $G_\delta$-neighborhood of $\chi$ in $K$, and these sets form a 
base for the $G_\delta$-topology on $K$.

If $\chi \in \bar H^\delta$, then $N(\chi,C)\cap H \neq \emptyset$ for 
every countable subset $C$ of $A$, and therefore for every $C$ there is
$\psi_C \in H$ such that $\psi_C \rest_C=\chi \rest_C$. In particular,
$\chi\rest_C$ is continuous with respect to the induced 
topology $\tau_H\rest_C$, and hence $\chi$ is countably continuous.

Conversely, suppose that $\chi$ is countably continuous, and let 
$C\subseteq A$ be countable. By replacing $C$ with the subgroup that it 
generates, we may assume that $C$ is a countable subgroup of $A$, because
$N(\chi,\langle C \rangle)= N(\chi,C)$.
One has $\chi\rest_C \in H \rest_C$, because
$\chi\rest_C$ is continuous with respect to $\tau_H\rest_C$
and $(C, \tau_H\rest_C)$ is a precompact.
Therefore,  $H \cap N(\chi,C) \neq \emptyset$, as desired. 
\end{proof}

We proceed with presenting a characterization of $\mathfrak l$-closed and 
$\mathfrak l$-dense
subgroups of compact Hausdorff abelian groups. To that
end, we first formulate explicitly an argument that was used by Comfort,
Hern\'andez and Trigos-Arrieta \cite{ComfHernTrig}.  In order to do so, 
the following terminology is 
needed: For a topological space $X$, a {\em zero-set in $X$} is a set of 
the form $Z(f)=f^{-1}(0)$, where $f$ is a real-valued continuous function 
on $X$. A subset  $Y \subseteq X$ is {\em $z$-embedded} in $X$ if for 
every zero-set $Z$ in $Y$ there is a zero-set $W$ in $X$ such that 
$Z=W\cap Y$; $Y$ is {\em $C$-embedded} in $X$ if every continuous 
real-valued function on $Y$ extends continuously to $X$.
One says that $X$ is an {\em $Oz$-space} if every open 
subset of $X$ is $z$-embedded.

\begin{proposition} \label{prop:HW-Gdelta}
Let $G$ be a Hausdorff group whose (Ra\u{\i}kov-)completion 
$\widetilde G$ is locally compact (in other words, $G$ is locally 
precompact), 
and let $H$ be its subgroup. Then: 

\begin{list}{{\rm (\alph{enumi})}}
{\usecounter{enumi}\setlength{\labelwidth}{25pt}\setlength{\topsep}{-10pt}
\setlength{\itemsep}{-4pt} \setlength{\leftmargin}{25pt}}

\item
$H$ is $z$-embedded in $\tilde G$;

\item
$H$ is $C$-embedded in $\bar H^\delta$; 

\item 
$\upsilon H \subseteq \upsilon G$ and
$\bar H^\delta= \upsilon H \cap G$, where
$\upsilon$ stands for the Hewitt-realcompactification;

\item
if $G$ is realcompact, then $\bar H^\delta=\upsilon H$.

\end{list}

\end{proposition}

Note that in contrast to the rest of this paper, in this proposition, the
group $G$ is not assumed to be abelian. We claim no novelty for the proof
of Proposition~\ref{prop:HW-Gdelta} (which heavily relies on results
of Blair and Hager), and it is included here for the sake of
completeness.

\begin{flemma}[{\cite[1.1(a)-(b)]{BlaHag2}, \cite[2.1, 5.1]{Blair}}]

\label{flemma:BlaHag}

\mbox{ }

\begin{list}{{\rm (\alph{enumi})}}
{\usecounter{enumi}\setlength{\labelwidth}{25pt}\setlength{\topsep}{-10pt}
\setlength{\itemsep}{-4pt} \setlength{\leftmargin}{25pt}}

\item
Every $z$-embedded $G_\delta$-dense subset of a topological space $X$ is 
$C$-embedded in $X$.

\item
If $Y$ is $z$-embedded in $X$, then $\upsilon Y \subseteq \upsilon X$
and $\upsilon Y$ is the $G_\delta$-closure of $Y$ in $\upsilon X$.

\item
A topological space $X$ is normal if and only if every closed subset of 
$X$ is $z$-embedded.

\item
A topological space $X$ is an $Oz$-space if and only if every dense subset 
of $X$ is $z$-embedded.

\end{list}

\end{flemma}

\begin{proof}
(a) The completion $\widetilde G$ is locally compact, 
and thus normal. By Lemma~\ref{flemma:BlaHag}(c), the closed
subgroup $\widetilde H$ of $\widetilde G$ is $z$-embedded in 
$\widetilde G$. By a result of Ross and Stromberg \cite{RossStro},
$\widetilde H$ is an $Oz$-space, because it is locally compact (as a 
closed subgroup of $\widetilde G$). Thus, the dense subgroup $H$ is 
$z$-embedded in $\tilde H$ (cf.  Lemma~\ref{flemma:BlaHag}(d)). Therefore,
$H$ is $z$-embedded in $\tilde G$.

(b) It follows from (a) that $H$ is $z$-embedded in $\bar H^\delta$, and 
therefore it is $C$-embedded in it, by  Lemma~\ref{flemma:BlaHag}(a).

(c) It follows from (a) that $H$ is $z$-embedded in $G$, and therefore the 
statement follows from Lemma~\ref{flemma:BlaHag}(b).

(d) follows from (c).
\end{proof}

\begin{corollary} \label{cor:l:misc}
Let $K\in \mathsf{Ab(CompH)}$ and $H\leq K$.

\begin{list}{{\rm (\alph{enumi})}}
{\usecounter{enumi}\setlength{\labelwidth}{25pt}\setlength{\topsep}{-10pt}
\setlength{\itemsep}{-4pt} \setlength{\leftmargin}{25pt}}

\item
If $H$ is closed, then it is $\mathfrak l$-closed in $K$.

\item
$H$ is $\mathfrak l$-closed in $K$ if and only if $H$ is realcompact.

\item
$H$ is $\mathfrak l$-dense in $K$ if and only if it is dense in $K$ and 
pseudocompact.

\item
$K$ is metrizable if and only if every subgroup of $K$ is 
$\mathfrak l$-closed.

\end{list}

\end{corollary}

\begin{proof}
In light of Theorem~\ref{thm:char:l}, (a) is obvious, 
(b) follows from Proposition~\ref{prop:HW-Gdelta}(d), and 
(c) is a reformulation of a famous theorem of Comfort and Ross
\cite[1.2]{ComfRoss2}.

(d) If $K$ is metrizable, then every subgroup of $K$ is realcompact,
and thus, by (c), $\mathfrak{l}$-closed. Conversely, if every subgroup of 
$K$ is $\mathfrak{l}$-closed, then $K$ is hereditarily realcompact (by 
(b)), and therefore metrizable (cf.~\cite[3.3]{HernMaca}).
\end{proof}

We conclude this section by observing that a result of Hern\'andez and
Macario can be obtained as an easy consequence of the foregoing
discussion:

\begin{fcorollary}[{\cite[3.2]{HernMaca}}] \label{fcor:real:ccc}
Let $K\in\mathsf{Ab(CompH)}$ and set $A=\hat K$.
A subgroup $H \leq K$ is realcompact if and only if every countably 
continuous character of $(A,\tau_H)$ is continuous.
\end{fcorollary}

\begin{proof}
As we explained in the proof of Theorem~\ref{thm:char:l}, $\mathfrak
l_K(H)$ is the set of countably continuous characters of $(A,\tau_H)$ (cf.  
Lemma~\ref{lemma:f}(b) and Lemma~\ref{lemma:l:cc}). On the other hand, by
Corollary~\ref{cor:l:misc}(b), $H$ is realcompact if and only if it is
$\mathfrak l$-closed in $K$.
\end{proof}


\section{Application II: The $\MakeLowercase{\mathfrak{g}}$-closure}

\label{sect:gclosure}

In this section, we apply the results of section~\ref{sect:cont}
to $\mathfrak{g}$, and  prove 
Theorems~\ref{thm:positive}-\ref{thm:gc:MAP}.
As shown in Example~\ref{example:S}, $\mathfrak{g}$ is the
closure operator associated with the continuity type $\mathcal{S}$
(assigning to a group the filter-bases generated by sequences).
Furthermore, $\mathcal{S}$-continuous homomorphisms are precisely
the sequentially continuous ones. 

\begin{theorem} \label{thm:g}
Let $K \in \mathsf{Ab(CompH)}$ and $A=\hat K$.
A subgroup $H \leq K$ is:

\begin{list}{{\rm (\alph{enumi})}}
{\usecounter{enumi}\setlength{\labelwidth}{25pt}\setlength{\topsep}{-10pt}
\setlength{\itemsep}{-4pt} \setlength{\leftmargin}{25pt}}

\item

$\mathfrak{g}$-closed in $K$ if and only if every 
sequentially continuous character of $(A,\tau_H)$ is continuous;

\item
$\mathfrak{g}$-dense in $K$ if and only if $(A,\tau_H)$ has no
non-trivial convergent sequences.

\end{list}
\end{theorem}

Note that part (b) of Theorem~\ref{thm:f} appeared earlier
in \cite[5.11]{BarDikMilWeb2}.

\begin{proof}
(a) is Theorem~\ref{thm:f} applied to $\mathfrak{g}$.

(b) If $(A,\tau_H)$ has no non-trivial convergent sequences, then
all its characters are sequentially continuous, 
and thus, by Lemma~\ref{lemma:f}(b), $\mathfrak g_K(H)=K$.
Conversely, suppose that $H$ is $\mathfrak{g}$-dense in $K$---in 
other words, $\mathfrak g_K(H)=K$. Then,
by Lemma~\ref{lemma:f}(c), the identity map
$(A,\tau_H)\rightarrow (A,\tau_{K})$ is sequentially continuous.
Thus, if $\{a_n\}$ is a convergent sequence with respect to
$\tau_H$, then it must converge in $\tau_K$. On the other hand,
one has $(A,\tau_K)=A^+ \subseteq bA$, and it is well known that 
$A^+$ has no convergent sequences, because $A$ is discrete abelian
(cf. \cite{Flor} and \cite{Glick}).
\end{proof}

It follows from the definition of $\mathcal{S}$ 
that it is regular, and every
dually embedded subgroup $N$ of a group $G\in\mathsf{Ab(Top)}$
is $\mathcal{S}$-dually embedded (because if
$\{u_n\}$ is a sequence in $\hat N$,
then there is a sequence $\{v_n\}$ in $\hat G$ such that
$\hat i_N(v_n)=u_n$). Therefore, the following well-known and basic properties of 
$\mathfrak{g}$ are immediate consequences of Proposition~\ref{prop:contype}
(cf.~\cite[2.9,2.11]{DikMilTon}):

\begin{proposition} \label{prop:gcl:basic} For every $G\in\mathsf{Ab(Top)}$:

\begin{list}{{\rm (\alph{enumi})}}
{\usecounter{enumi}\setlength{\labelwidth}{25pt}\setlength{\topsep}{-10pt}
\setlength{\itemsep}{-4pt} \setlength{\leftmargin}{25pt}}

\item
$\mathfrak{g}_G(\{0\})=\mathbf{n}(G)$;

\item
$\mathfrak{g}_G(H)=\mathfrak{g}_{G^+}(H)$ for every $H\leq G$;

\item
if $N\leq G$ is dually embedded in $G$, then 
$\mathfrak{g}_N(H)=\mathfrak{g}_G(H)\cap N$ for every $H \leq N$.
\qed
\end{list}
\end{proposition}

Now we are ready to prove the theorems that were formulated in
section~\ref{sect:intro}:

\begin{proof}[Proof of Theorem~\ref{thm:positive}.]
Let $H=\{h_1,h_2,\ldots\}$ be countable.  Since $K$
is compact, it follows from the Pontryagin
duality that its subgroups are dually embedded, and so
$\mathfrak g_{\bar H}(H)=\mathfrak g_K(H) \cap \bar H$ 
(by Proposition~\ref{prop:gcl:basic}). Every closed subgroup
of $K$ is $\mathfrak{g}$-closed (cf. \cite[2.12]{DikMilTon}),
and thus one has $g_K(H) \leq \bar H$. Therefore,
$\mathfrak g_{\bar H}(H)=\mathfrak g_K(H)$, and
without loss of generality, we may assume that $H$ is dense in $K$. 
In particular,  $(A,\tau_H)$ is Hausdorff.  
The topology $\tau_H$ is the initial one with respect to the homomorphism 
$(h_i)\colon A \rightarrow\mathbb{T}^\omega$. 
Thus, $(A,\tau_H)$ is metrizable, and in particular, every
sequentially continuous map on $(A,\tau_H)$ is continuous.
Therefore, by Theorem~\ref{thm:g}(a), $H$ is $\mathfrak{g}$-closed.
\end{proof}

\begin{proof}[Proof of Theorem~\ref{thm:gc-realcompact}.]
Let $H$ be a $\mathfrak g$-closed subgroup of $K$.
Since  $\mathcal{S}(\hat K) \subseteq \mathcal{C}(\hat K)$
(where $\mathcal{S}$ is as in Example~\ref{example:S}), one has
$H\leq \mathfrak l_K(H) \leq \mathfrak g_K(H) =H$, and thus
$H$ is $\mathfrak l$-closed in $K$. Therefore, by 
Corollary~\ref{cor:l:misc}(b), $H$ is realcompact. 
\end{proof}

\begin{proof}[Proof of Theorem~\ref{thm:gc:MAP}.]
(i) $\Rightarrow$ (ii): 
Let $H$ be a countable subgroup of $G$. 
By Theorem~\ref{thm:positive},  
$\rho_G(H)$ is $\mathfrak{g}$-closed in the Bohr-compactification $bG$, 
because it is countable---in other words,  
$\mathfrak g_{bG}(\rho_G(H))=\rho_G(H)$. Since
$G$ is maximally almost periodic (i.e., $\rho_G$ is injective),
one has $\rho_G^{-1}(\mathfrak g_{bG} (\rho_G(H)))=H$. Hence,
by Theorem~\ref{thm:reduce:bohr}, $H$ is $\mathfrak{g}$-closed
in $G$.


(ii) $\Rightarrow$ (iii) is trivial.

(iii) $\Rightarrow$ (i):
By (iii), the (cyclic) subgroup $\{0\}$ is $\mathfrak{g}$-closed
in $G$. Therefore, by Proposition~\ref{prop:gcl:basic}(a),
$\mathbf{n}(G)=\mathfrak g_G(\{0\})=\{0\}$, as desired.
\end{proof}

\section*{Acknowledgments}

I am deeply indebted to Professor Dikran Dikranjan, who introduced me 
to this topic and spent so much time corresponding with me about it.
Without his attention and encouragement, I would have struggled to 
carry out this research.

I wish to thank Professor Horst Herrlich, Professor Salvador
Hern\'andez and Dr. Gavin Seal for the valuable discussions 
and helpful suggestions that were
of great assistance in writing this paper.

\bibliography{notes,notes2,notes3}

\begin{samepage}

{\bigskip\bigskip\noindent 
Department of Mathematics and Statistics\\
Dalhousie University\\
Halifax, B3H 3J5, Nova Scotia\\
Canada

\nopagebreak
\bigskip\noindent{\em e-mail: lukacs@mathstat.dal.ca} }
\end{samepage}

\end{document}